\UseRawInputEncoding
\documentclass{amsart}
\usepackage{amsmath}
\usepackage{amssymb}
\newtheorem{thm}{Theorem}[section]
\newtheorem{cor}[thm]{Corollary}
\newtheorem{lem}[thm]{Lemma}
\newtheorem{prop}[thm]{Proposition}
\newtheorem{cons}[thm]{Construction}
\newtheorem{ex}[thm]{Example}
\theoremstyle{definition}
\newtheorem{defn}[thm]{Definition}
\newtheorem{notn}[thm]{Notation}
\theoremstyle{remark}
\newtheorem{rem}[thm]{Remark}
\long\def\Thm#1{\begin{thm} #1 \end{thm}}
\long\def\Cor#1{\begin{cor} #1 \end{cor}}
\long\def\Lem#1{\begin{lem} #1 \end{lem}}
\long\def\Prop#1{\begin{prop} #1 \end{prop}}

\long\def\Def#1{\begin{defn} #1 \end{defn}}

\long\def\Rem#1{\begin{rem} #1 \end{rem}}
\long\def\Ex#1{\begin{ex} #1 \end{ex}}

\def\Sect{\section}
\def\Rarr#1#2{\xrightarrow[#2]{#1}}

\def\Darr#1#2{{\scriptstyle #1}\downarrow{\scriptstyle #2}}

\long\def\Ref#1#2#3#4#5#6{
\bibitem{#1}
{\rm #2,}
\textit{#3.}
{\rm #4}
\textbf{#5}
{\rm #6.}
}
\long\def\Refb#1#2#3#4{
\bibitem{#1}
{\rm #2,}
\textit{#3.}
#4.
}
\def\Qq{{\mathbb Q}}
\def\Zz{{\mathbb Z}}
\def\Nn{{\mathbb N}}
\def\Rr{{\mathbb R}}
\def\Cc{{\mathbb C}}
\def\SS{\mathfrak{S}}

\def\LL{\mathcal{L}}

\def\into{\hookrightarrow}
\def\iso{\cong}
\def\leq{\leqslant}
\def\geq{\geqslant}
\def\hcf{{\rm hcf}}

\def\comp{\mathbin{\mathchoice
{\circ}
{{\scriptstyle\circ}}
{{\scriptscriptstyle\circ}}
{{\scriptscriptstyle\circ}}
}}
\def\st{\mid}
\def\Fix{{\rm Fix}}

\def\Hom{{\rm Hom}}
\def\map{{\rm map}}

\def\tr{{\rm tr}}

\def\phi{\varphi}

\def\Alb{{\rm Alb}}

\def\DD{\mathcal{D}}

\def\PP{\mathcal{P}}

\def\Bigwedge{\Lambda}
\def\SS{\mathfrak{S}}
\begin{document}

\title{Dold indices and symmetric powers}
\author{M.~C.~Crabb}
\address{%
Institute of Mathematics\\
University of Aberdeen \\
Aberdeen AB24 3UE, UK}
\email{m.crabb@abdn.ac.uk}
\begin{abstract}
Results of Macdonald and Dold from the 1960s and '70s
expressing the Lefschetz numbers of symmetric powers
of a self-map of a compact ENR in terms of the Lefschetz
numbers of iterates of the map are extended using
the notion of a Lefschetz-polynomial functor.
Configuration spaces and Borsuk-Ulam symmetric products,
as well as symmetric powers, are treated as examples of the
general method.
\end{abstract}
\date{March 2025}
\maketitle

\Sect{Introduction}
In 1962 Macdonald \cite{IM} showed that the rational cohomology
Poincar\'e polynomial of
the $k$th symmetric power $S^kM$ of a finite complex $M$
is determined by the Poincar\'e polynomial of $M$.
Ten years later Dold \cite{AD} proved that, for a continuous map
$f: M\to M$ inducing linear maps $A_i:H^i(M;\,\Qq )\to H^i(M;\,\Qq)$,
the characteristic rational function 
$\prod_i \det (1-tA_i))^{(-1)^i}\in\Qq (t)$ of $f$ 
determines the characteristic rational function of 
the $k$th symmetric power of $f$, $S^kf : S^kM\to S^kM$.
In terms of Lefschetz numbers, Dold's result says that the
Lefschetz numbers of the powers of $S^kf$ are determined
by the Lefschetz numbers $L(f^k)$ of the iterates $f^k:M\to M$ of $f$. 
The Dold indices $\DD^m(f)\in\Zz$,
$m\geq 1$, introduced in \cite{AD2}, count the fixed points of the 
iterates of $f$ and are related to the Lefschetz numbers
through the associated zeta-function defined as
$$
Z(f; q)=\prod_{m\geq 1}(1-q^m)^{\DD^m(f)}\in\Qq [[q]]
$$
by the formal power series identity
$$
Z(f; q)
=\exp\big(-\sum_{k\geq 1}\frac{L(f^k)}{k}q^k\big)\in\Qq [[q]].
$$
If, for some $k\geq 1$, the Lefschetz index $L(S^kf)$
is non-zero, then there exists some non-empty 
finite subset $A\subseteq M$ of cardinality at most $k$
such that $f(A)=A$.

The results of Macdonald and Dold, reformulated in the 
following theorem, provide the starting point for this note. 
A short direct proof is given in Section 3.
The $0$th power $S^0M$ is understood to be a point.
\Thm{\label{MD}
{\rm (Macdonald \cite{IM}, Dold \cite[Theorem 5.9]{AD}
and \cite[Theorem 4.1]{AD3}).}
Let $f: M\to M$ be a (continuous) self-map of a compact ENR $M$.
Then the Lefschetz indices of the maps $S^kf :S^kM \to S^kM$ 
are given by
$$
\sum_{k\geq 0}  L(S^kf : S^kM\to S^kM )\, q^k
= Z(f; q)^{-1}\in \Qq [[q]].
$$
}
In Section 2 we introduce the notion of a `Lefschetz-polynomial
functor',
which allows us to reduce the verification of identities
of this type to the 
combinatorial special case of a self-map of a finite set.
The proof that certain natural functors are Lefschetz-polynomial
is given in Section 5.  Symmetric powers, configuration spaces 
and Borsuk-Ulam symmetric products are treated as examples
of the general method. 
\Sect{Lefschetz-polynomial functors}
Motivated by Dold's work \cite{AD, AD1} on non-additive functors
(see also \cite{vall, dB}),
we introduce the idea of a {\it Lefschetz-polynomial functor}
tailored to our specific investigation of Lefschetz numbers.

Consider a pointed compact ENR (Euclidean Neighbourhood Retract) $X$
and a pointed self-map $f : X\to X$. 
The {\it reduced Lefschetz number} of $f$
is defined in terms of the Lefschetz number
$L(f)$ as $\tilde L(f)=L(f)-1\in\Zz$.
(Its significance for the
existence of fixed points is that, if $f$ maps a neighbourhood
of the basepoint $*$ to the basepoint and $\tilde L(f)$ is non-zero,
then $f$ has a fixed-point other than $*$.) 
Similarly, we introduce the {\it reduced Dold indices}
$\tilde\DD_1(f)=\DD_1(f)-1$ and $\tilde \DD_m(f)=\DD_m(f)$
for $m>1$, so that the {\it reduced zeta function}
$$
\tilde Z(f; q)=\prod_{m\geq 1}(1-q^m)^{\tilde\DD^m(f)}\in\Qq [[q]]
$$
is equal to
$$
\exp\big(-\sum_{k\geq 1}\frac{\tilde L(f^k)}{k}q^k\big)\in\Qq [[q]].
$$
\Def{A functor $\Phi$ from the category of pointed compact ENRs and
basepoint-preserving
continuous maps to itself is said to be {\it Lefschetz-polynomial
of degree at most $k$}
if there is a polynomial $L_\Phi (t_1,\ldots ,t_k)\in
\Qq [t_1,\ldots ,t_k]$ of weighted degree at most $k$, where
$t_i$ has weighted degree $i$, such that for any pointed
self-map $f: X\to X$
$$
\tilde L(\Phi f :\Phi X\to\Phi X) =
L_\Phi (\tilde\DD_1(f),\ldots ,\tilde\DD_k(f))\in\Zz\, .
$$
}
\Rem{If $\Phi$ is Lefschetz polynomial of degree at most $k$, then
it is Lefschetz polynomial of degree at most $l$ for any $l\geq k$.
We shall write, without risk of ambiguity,
$L_\Phi (t_1,\ldots ,t_l)=L_\Phi (t_1,\ldots ,t_k)$.
}
\Lem{The polynomial $L_\Phi$ is uniquely determined by $\Phi$.
It is numerical, in the sense that $L_\Phi (d_1,\ldots ,d_k)\in\Zz$
if $d_i\in\Zz$, and it is determined by its values on
$(\tilde\DD_m(f))$ for any self-map $f$ of a finite pointed set $X$.
}
\begin{proof}
For any non-negative integers $d_1,\ldots ,d_k$ there is a finite set
$M$ and a self-map $f :M\to M$ with exactly $d_i$ orbits of length $i$
for each $i=1,\ldots ,k$ and no orbits of greater length. Take $X=M_+$
to be the pointed space obtained by adjoining a basepoint to $M$ and
consider the induced pointed self-map of $X$.
The values $L_\Phi (d_1,\ldots ,d_k)$ for $d_i\geq 0$ determine the
polynomial $L_\Phi (t_1,\ldots ,t_k)$.
(It is enough to observe that, if a polynomial $p(t)\in\Qq [t]$
takes integer values $p(d)\in\Zz$ for non-negative integers
$d$, it is a $\Zz$-linear combination of binomial coefficients
$\binom{t}{r}$, $r\geq 0$.)
\end{proof}
\Lem{\label{add}
{\rm (Additivity).}
Suppose that $\Phi',\, \Phi$ and $\Phi''$ are three functors
related by natural transformations $\Phi'\to\Phi\to\Phi''$
such that for any pointed compact ENR $X$
$$
\Phi' (X) \to \Phi (X)\to \Phi'' (X)
$$
is a homotopy cofibre sequence. If any two of the functors
$\Phi',\, \Phi$ and $\Phi''$
are Lefschetz-polynomial of degree at most $k$, then
so, too, is the third and
$$
L_\Phi (t_1,\ldots ,t_k)=L_{\Phi'}(t_1,\ldots ,t_k)
+L_{\Phi''}(t_1,\ldots ,t_k).
$$
}
\begin{proof}
For a map $f:X\to X$ we have a commutative ladder of long exact 
sequences
$$
\begin{matrix}
\cdots\to&\tilde H^i(\Phi''(X);\,\Qq )&\to&\tilde H^i(\Phi (X);\,\Qq )
&\to&H^i(\Phi'(X);\,\Qq )&\to\cdots\\
&\Darr{(\Phi'f)^*}{}&&\Darr{(\Phi f)^*}{}&&\Darr{(\Phi''f)^*}{}\\
\cdots\to&\tilde H^i(\Phi''(X);\,\Qq )&\to&\tilde H^i(\Phi (X);\,\Qq )
&\to&H^i(\Phi'(X);\,\Qq )&\to\cdots
\end{matrix}
$$
Hence, $\tilde L(\Phi f)=\tilde L(\Phi'f)+\tilde L(\Phi''f)$.
\end{proof}
\Lem{\label{mult}
{\rm (Multiplicativity).}
Suppose that $\Phi_1$ and $\Phi_2$ are 
Lefschetz-poly\-nomial functors
of degree at most $k_1$ and $k_2$ respectively.
Then the smash product $\Phi (-)=\Phi_1(-)\wedge\Phi_2(-)$
is Lefschetz-polynomial of degree at most $k=k_1+k_2$ and
$$
L_\Phi (t_1,\ldots ,t_k)=L_{\Phi_1}(t_1,\ldots ,t_{k_1})
\cdot L_{\Phi_2}(t_1,\ldots ,t_{k_2}).
$$
}
\begin{proof}
For a self-map $f:X\to X$, we have $\tilde L(\Phi f)=
\tilde L(\Phi_1f)\cdot \tilde L(\Phi_2f)$.
\end{proof}
\Ex{Suppose that $A_0,\ldots ,\, A_k$ are pointed compact
ENRs. Then the functor $\Phi X =A_0\vee(A_1\wedge X)\vee\cdots
\vee (A_i\wedge \Bigwedge^iX)\vee\cdots \vee (A_k\wedge\Bigwedge^kX)$
is polynomial of degree at most $k$ and
$L_\Phi (t_1,\ldots ,t_k)=a_0+a_1t_1+\cdots +a_kt_1^k$,
where $a_i\in\Zz$ is the Euler characteristic of $A_i$.
}
\Lem{\label{iter}
Suppose that a functor $\Phi$ is Lefschetz-polynomial of 
degree at most $k$. Then there exist unique numerical
polynomials 
$$
D_m^\Phi (t_1,\ldots ,t_{mk})\in\Qq [t_1,\ldots ,t_{mk}]
$$
of weighted degree at most $mk$ such that for any pointed
self-map $f: X\to X$
$$
\tilde \DD_m(\Phi f :\Phi X\to\Phi X) =
D_m^\Phi (\tilde\DD_1(f),\ldots ,\tilde\DD_{mk}(f))\in\Zz\, .
$$
The polynomial $D_m^\Phi$ is determined by its values on
$(\tilde\DD_n(f))$ for any self-map $f$ of a finite pointed set $X$,
}
\begin{proof}
For a map $f: X\to X$, the Dold index $\tilde\DD_m(\Phi f)$
is a rational linear combination of the Lefschetz numbers
$\tilde L ((\Phi f)^j)$, $j\, |\, m$.
Since $(\Phi f)^j=\Phi (f^j)$, $\tilde\DD_m(\Phi f)$
is a linear combination of $L_\Phi (\tilde\DD_1(f^j),\ldots ,
\tilde\DD_k(f^j))$, $j\, |\, m$.
Each $\tilde\DD_i(f^j)$ is an integral linear combination
$$
\tilde\DD_i(f^j)=\sum_{i\cdot\hcf (j,n)=n} \tilde\DD_n(f).
$$
Thus $\tilde\DD_m(\Phi f)$ is expressed as a polynomial 
of (weighted)
degree at most $mk$ in the $\tilde\DD_n(f)$.
\end{proof}
\Lem{\label{comp}
{\rm (Composition).}
Suppose that $\Phi$ and $\Psi$ are Lefschetz-polynomial functors
of degree at most $k$ and $l$ respectively. Then the
composition $\Psi\comp\Phi$ is Lefschetz-polynomial of degree
at most $k^l$.
}
\begin{proof}
This follows directly from Lemma \ref{iter},
because 
$$
\tilde L (\Psi\comp\Phi (f))
=L_\Psi (\DD_1(\Phi (f)),\ldots ,\DD_l(\Phi (f)))
$$
and each $\DD_i(\Phi (f))$ for $i\leq l$
is a polynomial in $\DD_j(f)$ for $j\leq kl$.
\end{proof}
\Ex{The simplest example of a functor $\Phi$ which is not 
Lefschetz-polynomial is the pointed set of components:
$\Phi X =\pi_0(X)$.

Another example is the Albanese torus 
$\Alb (X)=(\Rr\otimes A(X))/A(X)$, 
where $A(X)$ is the finitely generated free 
abelian group $\Hom (\tilde H^1(X;\, \Zz ),\Zz )$,
that is, the abelianization of the fundamental group $\pi_1(X)$.
}
\Sect{Symmetric powers}
Let $M$ be a compact ENR. For an integer $k\geq 1$ we write
$S^kM =M^k/\SS_k$ for the $k$-th {\it symmetric power},
the orbit space of the action of the symmetric group $\SS_k$
on $M^k$.
We can think of a point of $S^kM$ as a formal sum
$\sum_{x\in M} m_xx$, where $m_x\in\Nn$, $m_x=0$ for all
but finitely many $x\in M$ and $\sum m_x=k$.

\begin{proof}[Proof of Theorem \ref{MD}]
We give an elementary proof following Macdonald's
argument in \cite{IM}
to obtain a result about the Poincar\'e polynomial
generating function
$$
P(f; q,T)=
\sum_{k\geq 0} q^k\sum_{i\geq 0} (-T)^i \tr \{ S^kf:
H^i(S^kM;\,\Qq )\to H^i(S^kM ;\,\Qq )\},
$$
in the ring $(\Qq [T])[[q]]$ of formal power series over the
polynomial ring $\Qq [T]$,
that strengthens Dold's Theorem 5.9 in \cite{AD}. 

Let us write $V_i=H^i(M;\,\Qq )$ and $A_i=f^*: V_i\to V_i$.
The pullback $H^*(S^kM;\,\Qq )$ $\to H^*(M^k;\,\Qq)$
maps isomorphically onto the invariant subalgebra\footnote{Macdonald
refers to Grothendieck's T\^ohoku paper \cite{Grothendieck},
Theorem 5.3.1 and the Corollary to Proposition 5.2.3.
A proof in the spirit of the present paper
is outlined in the Appendix.}
$H^*(M^k;\, \Qq )^{\SS_k}$. So $P(f; q,T)$ can be expressed in
terms of $A=(A_i)$ and $V=(V_i)$ as
$$
P(A; q,T)=
\sum_{k} q^k\sum_i (-T)^i
\textstyle{\tr \{ (\bigotimes^k A)_i:
(\bigotimes^k V)_i^{\SS_k}\to (\bigotimes^k V)_i^{\SS_k})\}.}
$$

This definition makes sense for any endomorphism $A$ of a 
finite-dimensional $\Zz$-graded $\Qq$-vector space $V$ with
$V_i=0$ for $i<0$. It is straightforward
to check that for another endomorphism $A' :V'\to V'$
of the same type,
$$
P(A\oplus A'; q,T) =P(A; q,T)\cdot P(A'; q,T),
$$
and that if, for some $j$, $V_i=0$ for $i\not=j$,
then, if $j$ is even,
$$
P(A; q,T) = \det (1-qT^jA_j)^{-1} =\det (1-qT^jA_j)^{-(-1)^j},
$$
and, if $j$ is odd,
$$
P(A; q,T)=\det (1-qT^jA_j)
=\det (1-qT^jA_j)^{-(-1)^j}.
$$
(It is enough to check the identity when $A_j$ is a diagonal
matrix.)
We conclude that, in general,
$$
P(A; q,T)=\prod_j \det (1-qT^jA_j)^{-(-1)^j}\, .
$$

The zeta function $Z(f; q)$ is equal to 
$$
\exp (-\sum_{k\geq 1}\frac {q^k}{k} \sum_j
(-1)^j \tr (A_j^k))
=\prod_j \det (1-qA_j)^{(-1)^j}\, .
$$
So, finally, $Z(f; q)=P(f; q,1)^{-1}$.
\end{proof}
A proof from an axiomatic point of view can be found in \cite{gomez}.

\smallskip

For an integer $l\geq 0$, let $C_l^k(M)\subseteq S^kM$
be the subspace of points $\sum m_xx$ such that $m_x\leq l$
for all $x\in M$. It is an open subspace. In particular,
$C_1^k(M)$ is the configuration space, $C^k(M)$, of $k$-element subsets
of $M$, and $C_0^k(M)=\emptyset$. 
The complement of $C_l^k(M)$ in $S^kM$ is a closed sub-ENR
(see Remark \ref{enr}),
and the quotient $S^k_lM=S^kM/(S^kM-C_l^k(M))$ is a pointed
compact ENR, which may be identified with the one-point
compactification of $C_l^k(M)$. In particular, for $l\geq k$,
$S_l^kM=(S^kM)_+$ is the pointed space obtained by
adjoining a disjoint basepoint to $S^kM$, and, if $k\geq 1$,
$S_0^kM=*$. When $k=0$, we take $S_l^0=S^0$ for all $l$.

A self-map $f :M\to M$ induces pointed maps 
$$
\begin{matrix}
(S^kM)_+&=S^k_\infty M&\to\cdots\to&
S^k_lM & \to\cdots\to & S^k_1M &\to &S^k_0M &=* \\
\Darr{(S^kf)_+}{}&&&
\!\!\!\!\!\Darr{S^k_kf}{\phantom{S^k_kf}}\!\!\!\!\!&&
\!\!\!\!\!\Darr{S^k_1f}{\phantom{S^k_lf}}\!\!\!\!\!&&
\!\!\!\!\!\Darr{S^k_0f}{\phantom{S^k_0f}}\!\!\!\!\!&\\
(S^kM)_+&=S^k_\infty M&\to\cdots\to&
S^k_lM & \to\cdots\to & S^k_1M &\to &S^k_0M& =* 
\end{matrix}
$$
\Thm{\label{main}
Let $f :M\to M$ be a self-map of a compact ENR.
Then, for $l\geq 0$,
$$
\sum_{k\geq 0} \tilde L(S^k_lf : S^k_lM\to S^k_lM )\, q^k
= Z(f; q^{l+1})Z(f; q)^{-1}.
$$
}
\begin{proof}
We give the proof first when $M$ is a finite set. 
Then $f$ has $\DD^m(f)$ orbits of length $m$. The Lefschetz numbers
count the fixed points. Thus $\tilde L (S^k_lf)$ is the
number of points $\sum m_xx\in C^k_l(M)$ such that
$$
\sum m_xx =\sum m_x f(x).
$$
Such a fixed point can be written as a sum of orbits of $f$:
$$
\sum m_CC,
$$
where $0\leq m_C\leq l$ and $\sum m_C\# C =k$.
Thus
$$
\sum_k \tilde L (S^k_lf)q^k
=\prod_C (1+q^{\# C}+q^{2\# C}+\ldots +q^{l\# C})
$$
$$
=\prod_m (1+q^m+q^{2m}+\ldots +q^{lm})^{\DD_m (f)}
=\prod_m (1-q^{(l+1)m})^{\DD_m(f)}(1-q^m)^{-\DD_m(f)},
$$
which is $Z(f; q^{l+1})Z(f; q)^{-1}$.

The proof will be completed by extending the functor $S^k_l$
to pointed spaces and showing that it is Lefschetz-polynomial
of degree $k$. This is done in Section \ref{parti} (Lemma \ref{kl}).
\end{proof}
\Ex{\label{circle}
For an elementary, explicit, example, take $M$ to be a 
sphere of odd dimension and let 
$f: M\to M$ be a map of degree $d$.
Then
$$
\tilde L (S^k_lf)=\begin{cases}
1&\text{if $k=0$,}\\
(1-d)d^j&\text{if $j(l+1)<k< (j+1)(l+1)$,}\\
0&\text{if $l+1$ divides $k$.}
\end{cases}
$$
}
\begin{proof}
We have $Z(f; q) =(1-q)(1-dq)^{-1}$.
The Lefschetz numbers are expressed by Theorem \ref{main}
in terms of 
$$
Z(f; q^{l+1})Z(f; q)^{-1}= (1+q+\ldots +q^l)
(1-dq)(1-dq^{l+1})^{-1}
$$
$$
=1+ (1-d)(q+\ldots +q^l)(1-dq^{l+1})^{-1}\, .
$$
\vglue-\baselineskip
\end{proof}
\Rem{\label{cycle}
By pushing a configuration of $k$ particles on the circle
$M=\Rr /\Zz$ until the distances between adjacent particles are equal
(as described, for example, in \cite[II, Example 14.4]{FHT})
one can see that $C^k(M)$ is homeomorphic to the total space
$(\Rr /\Zz )\times_{\frac{1}{k}\Zz /\Zz } V_k$, where
$V_k$ is the reduced regular representation of the cyclic group
$\frac{1}{k}\Zz /\Zz$, of a
$(k-1)$-dimensional real vector bundle $\xi_k$ over the circle
$\Rr /\frac{1}{k}\Zz$. The bundle $\xi_k$ is
orientable if $k$ is odd and so trivial, non-orientable if
$k$ is even and so isomorphic to the direct sum of the Hopf line bundle
and a trivial bundle of dimension $k-2$.
Hence, $\tilde H^*(C^k(M);\,\Qq )=0$ if $k$ is even;
this is consistent with the vanishing of $\tilde L(S^k_1f)$
in Example \ref{circle}.
}
Both Macdonald and Dold generalized the symmetric
product $S^kM$ to the space $\map (K,M)/G$, where
$G$ is a finite group and $K$ is a finite $G$-set.
\Thm{\label{Gsymm}
Suppose that $G$ is a finite group and $K$ a finite
$G$-set. Then the functor $\Phi$ taking a pointed compact ENR $X$
to the one-point compactification of the orbit space
$\map (K,X-\{ *\})/G$ is Lefschetz-polynomial of degree 
at most $k=\# K$
and
$$ 
L_\Phi (t_1,\ldots ,t_k)=
\frac{1}{\# G}\sum_{g\in G}\left( \prod_{n\, |\, \# G}
(\sum_{m\, |\, n} mt_m)^{d_g(n)}\right),
$$ 
where $d_g(n)$ is the number of $g$-orbits of length $n$
in $K$, so that $k=\sum n \, d_g(n)$.
}
\begin{proof}
For a pointed compact ENR $X$, write 
$$
\map_0(K,X)=
\{a\in\map (K,X)\st*\in a(K)\}\subseteq \map (K,X). 
$$
It is a closed sub-$G$-ENR (see Remark \ref{enr}), 
and the quotient $\Phi X$
of $\map (K,X)/G$ by the subspace $\map_0(K,X)/G$
is a compact ENR. 
So, if $X=M_+$ is obtained by adjoining
a basepoint to a compact ENR $M$, we have
$\Phi (M_+) =(\map (K,M)/G)_+$.

For a finite set $M$ and map $f: M\to M$, we can count the fixed points
of $\Phi (f_+): (\map (K,M)/G)_+ \to (\map (K,M)/G)_+$.
This is equal to
$$
\frac{1}{\# G}\sum_{g\in G} \#\Fix (fg^{-1} :\map (K,M)
\to\map (K,M)).
$$
For a given $g\in G$,
we decompose $K$ as a union of $g$-orbits, $d_g(n)$ of length $n$,
and count the number of fixed points of $fg^{-1}$ in $\map (C,M)$ 
for an orbit $C$
of length $n$ as $\sum_{m\, |\, n} m\DD_m(f)$.

To show that $\Phi$ is Lefschetz-polynomial
(and confirm the formula already obtained)
we can follow Macdonald's argument using the expression
for $\tilde L(\Phi (f))$ as the average
$$
\frac{1}{\# G} \sum_{g\in G} \tilde L( f_*g^{-1} :
\map (K,X)/\map_0(K,X) \to\map (K,X)/\map_0(K,X))
$$
in $\Qq$-cohomology.
Taking $K$ to be the standard set $\{ 1,\ldots ,k\}$,
we recognize the quotient 
$\map (K,X)/\map_0(K,X)$ as the $k$-fold smash
product $\Bigwedge^k X$. Thus
$$
\tilde L(\Phi f)=
\frac{1}{\#G} \sum_g (-1)^i \tr ((\Bigwedge^kf^*)g^{-1}:
\tilde H^i(\Bigwedge^k X;\,\Qq )\to
 \tilde H^i(\Bigwedge^k X;\,\Qq )).
$$
Write $A_i=f^* : \tilde H^i(X;\,\Qq )\to \tilde H^i(X;\, \Qq )$.
Then the contribution from $g\in G$ is equal to the product
over the $g$-orbits $C$ of length $n_C$ of
$$
\sum_i (-1)^i \tr (A_i^{n_C})=\tilde L (f^{n_C})
=\sum_{m\, |\, n_C} m\tilde\DD_m(f)
$$
divided by $\# G$.
\end{proof}
\Sect{Borsuk-Ulam symmetric products}
In \cite{BU} Borsuk and Ulam introduced a different notion of the
$k$-fold {\it symmetric product} of a compact ENR $M$ as the space, which we shall write\footnote{The notation  ${\rm Sub}(M,k)$ is used in
\cite{handel} and $\exp_kM$ in \cite{tuff}.} as $P^kM$, of finite non-empty subsets of $M$ with at most
$k$ elements, topologized as the quotient
$$
P^kM=M^k/\sim,
$$
where $(x_1,\ldots ,x_k)\sim (y_1,\ldots ,y_k)$ if
and only if $\{ x_1,\ldots ,x_k\}=\{ y_1,\ldots ,y_k\}$.
(Thus, a subset $W\subseteq P^kM$ is a neigbourhood
of a point $A\in P^kM$ if and only if there exist
open sets $(U_x)_{x\in A}$ in $M$ with $x\in U_x$
such that every $B\in P^kM$ with $U_x\cap B\not=\emptyset$
lies in $W$.)
For $k=0$, $P^0M$ is empty.

The space $P^kM$ is a compact ENR, and there is an obvious inclusion
$P^{k-1}M\into P^kM$ for $k\geq 1$, with complement
$P^kM-P^{k-1}M=C^k(M)$, so that the quotient
$P^kM/P^{k-1}M$ is identified with $S^k_1M$.
(For the general topology, see \cite{jan, handel}.
It is enough to show that $P^kM$
is compact Hausdorff, locally contractible, and can be embedded as a subspace of some Euclidean space.)
\Prop{\label{prod}
Let $f: M\to M$ be a self-map of a compact ENR.
Then
$$
\sum_{k\geq 1} L(P^kf : P^kM\to P^kM)\, q^k = 
(1-q)^{-1}(Z(f; q^2)Z(f; q)^{-1}-1)\, .
$$
}
\begin{proof}
This is easily deduced from Theorem \ref{main} using the cofibre
sequence
$$
\begin{matrix}
P^{k-1}M &\Rarr{}{} & P^k M & \Rarr{}{} & S^k_1M \\
\noalign{\smallskip}
\Darr{P^{k-1}f}{}&&\Darr{P^kf}{}&& \Darr{S^k_1f}{}\\
\noalign{\smallskip}
P^{k-1}M &\Rarr{}{} & P^k M & \Rarr{}{} & S^k_1M 
\end{matrix}
$$
for $k\geq 1$.
\end{proof}
\Rem{When $M$ is a finite set and $f$ has $\DD_m(f)$ orbits of length
$m$, $L(P^kf)$ is equal to the sum of the coefficients of $q^j$ in
$\prod_m(1+q^m)^{\DD_m(f)}$ for $1\leq j\leq k$.
For a pointed compact ENR $X$ one can define a symmetric $k$-fold
smash product as the quotient of
the space $P^kX$ of non-empty subsets of $X$ of cardinality
at most $k$ by the closed subspace $P^k_0X$ consisting of 
those subsets that contain the basepoint $*$. 
It follows from the results in
Section \ref{parti} that this functor is Lefschetz-polynomial
of degree $k$.
Indeed, for $k>1$ the cofibre of the inclusion
$$
P^{k-1}X/P^{k-1}_0X \into P^kX/P^k_0X 
$$
can be identified with the one-point compactification
of the configuration space $C^k(X-\{ *\})$ of $k$-element
subsets of $X$ that do not contain $*$.
}
\Ex{Taking $f$ to be the identity, so that $Z(f; q)=(1-q)^\chi$,
where $\chi$ is the Euler characteristic of $M$, we see that
the Euler characteristic of $P^kM$ is equal to
$\sum_{1\leq j\leq k} \binom{\chi}{j}$.
{\rm (See \cite[Proposition 8]{sal} for a special case.)}
}
\Ex{\label{circle2}
We continue Example \ref{circle} in which $M$ is an 
odd-dimensional sphere.
It follows from Proposition \ref{prod} 
that $L(P^{2l-1}f)=L(P^{2l}f)=1-d^l$ for $l\geq 1$.
If $X$ is a pointed sphere of odd dimension,
a basepoint-preserving map $f:  X\to X$ induces a map
$P^kX/P^k_0X\to P^kX/P^k_0X$ with
reduced Lefschetz number equal to $-d^l$ if $k=2l-1$,
$0$ if $k=2l$.
}
\Rem{
When $X$ is the circle $\Rr /\Zz$ with $0$ as basepoint, 
the pointed space $P^kX/P^k_0X$ is homotopy
equivalent to $S^{2l-1}$ if $k=2l-1$ is odd, and is
contractible if $k$ is even.
Similarly, the pointed space $P^k_0X$ is contractible if
$k$ is odd and is homotopy equivalent to $S^{2l-1}$ if
$k=2l$, \cite[Theorem 2]{tuff2}.

(Here is a sketch of the argument. Let us write the standard
$k$-simplex as $\Delta^k=\{ (t_0,\ldots ,t_k)\in\Rr^{k+1}\st
t_i\geq 0,\, \sum_{i=0}^kt_i=1\}$. We have maps
$$
\Delta^k \Rarr{\alpha_k}{}P^kX/P^k_0X 
\Rarr{\beta_k}{} P^kX/(P^k_0X\cup P^{k-1}X) \Rarr{\gamma_k}{\iso} \Delta^k/
\partial \Delta^k
$$
determined by $\alpha_k(t_0,\ldots ,t_k)=[\{ x_1,\ldots ,x_k\}]$,
where $x_i =\sum_{0\leq j<i}t_j + \Zz$,
and the requirement that the composition be the quotient map
$\Delta^k\to\Delta^k/\partial\Delta^k$.
The map $\alpha_k$ restricts on the boundary to a map
$$
\partial \Delta^k \Rarr{\alpha_k |}{}
 P^{k-1}X/P_0^{k-1}X \Rarr{\gamma_{k-1}\comp
\beta_{k-1}}{} \Delta^{k-1}/\partial\Delta^{k-1}
$$
of degree $\sum_{j=0}^k (-1)^j$, that is, $1$ if $k$ is even, 
$0$ if $k$ is odd.

The description of $P^kX/P^k_0X$
follows by induction on $k$. If $k$ is odd and 
$P^{k-1}X/P_0^{k-1}$ is contractible, then $\gamma_k\comp\beta_k$
must be a homotopy equivalence. If $k$ is even and 
$\gamma_{k-1}\comp\beta_{k-1}$ is a homotopy equivalence,
then $\alpha_k$ has to be a homotopy equivalence.

In the same way for $P^k_0X$
we have maps 
$$
\Delta^{k-1}\Rarr{\alpha'_k}{} P^k_0X
\to P^k_0X/P^{k-1}_0X \to \Delta^{k-1}/\partial \Delta^{k-1},
$$
with
$\alpha'_k(t_0,\ldots ,t_{k-1})=\{ *,\, x_1,\ldots ,x_{k-1}\}$.
The restriction $\partial\Delta^{k-1} \to P^{k-1}_0X
\to \Delta^{k-2}/\partial\Delta^{k-2}$ has degree
$\sum_{j=0}^{k-1}(-1)^j$.)

It follows at once 
that $H^*(P^{2l-1}X;\,\Zz )\iso H^*(S^{2l-1};\, \Zz )
\iso H^*(P^{2l}X;\,\Zz )$, as was established in
\cite[Theorem 4]{tuff} (and in Bott's letter \cite{bott} to
Borsuk for $P^4X$). 
So Example \ref{circle2} reproves \cite[Theorem 7]{tuff}.
}
\Sect{\label{parti}
Partitions}
Let $G$ be a finite group and $K$ a finite $G$-set with $\# K=k$.
We write $\PP (K)$ for the set of partitions ordered by
refinement (so that the partition into $k$ $1$-element subsets
is the least element $0$ and the trivial partition $\{ K\}$ is the
greatest element $1$).
Consider a non-empty subset $\LL \subseteq\PP (K)$ closed under the
action of $G$ and taking refinement, that is, if $\pi\in\LL$,
then $g\pi\in\LL$ for all $g\in G$, and, if $\pi'\leq\pi$,
then $\pi'\in\LL$.

Any element $a\in \map (K,X)$ defines a partition $\pi (a)$ into
the non-empty subsets $a^{-1}(x)$, $x\in X$.
Define ${}_G\Phi^K_\LL X$ to be the quotient
$$
(\map (K,X)/(\{ a\st \pi (a)\notin\LL \}\cup \map_0(K,X)))/G\, .
$$

Equivalently, ${}_G\Phi^K_\LL X$ is the one-point compactification
of the configuration space
$$
{}_GC_\LL^K (X-\{ *\}) =\{ a\in\map (K,X-\{ *\})\st \pi (a)\in 
\LL \}/G\, .
$$
\Rem{\label{enr}
Using standard properties of ENRs as recalled below,
it is straightforward to check that
the orbit space ${}_G\Phi^K_\LL X$ is a compact ENR.

(a) Let $G$ be a finite group. If $A$ is a $G$-space such that
for each subgroup $H\leq G$ the fixed subspace $A^H$ is a compact
ENR, then $A$ is a $G$-ENR and the orbit space $A/G$ is a compact
ENR.
(b) Suppose that $(A_i)_{i\in I}$ is a finite family of closed
subspaces of a compact ENR $A$ such that each intersection
$\bigcap_{i\in J}A_i$  for $J\subseteq I$ is an ENR. 
Then the union $\bigcup_{i\in I}A_i$ is a compact sub-ENR of $A$.
(c) If $B$ is a closed sub-ENR of a compact ENR $A$, then
the quotient $A/B$ is a compact ENR.

For a partition $\lambda\in\PP (K)$, write
$X_\lambda = \{ a\in\map (K,X)\st \lambda \leq \pi (a)\}$
for the set of maps $a$ that are constant on each block
of $\lambda$. 
The verification that ${}_G\Phi^K_\LL X$ is a compact ENR
reduces to the observations:
(i) for a subgroup $H\leq G$, $\map (K,X)^K=\map (K/H,X)$;  
(ii) for $\lambda ,\, \mu \in \PP (K)$,
$X_\lambda \cap X_\mu = X_{\lambda \vee \mu}$. 
}
\Lem{\label{kl}
The functor ${}_G\Phi_\LL^K$ defined above is Lefschetz-polynomial
of degree at most $k$.
}
The proof given below reduces when $G$ is trivial to
the argument used in \cite[Lemma 2.5]{hainaut}.
\begin{proof}
We have established the result for the case
$\LL =\PP (K)$ in Theorem \ref{Gsymm}.
The proof will be completed by induction on $k=\# K$.

Suppose that $\LL\not=\PP (K)$.
Choose a minimal element $\lambda$ in $\PP (K)-\LL$.
Thus, $\lambda$ is not in $\LL$
but every proper refinement of $\lambda$ is in $\LL$.
Set $\LL^+=\LL \sqcup \{ g\lambda\st g\in G\}$.
We have a cofibre sequence
$$
\Phi'X \to {}_G\Phi^K_{\LL^+} X \to {}_G\Phi^K_\LL X
$$
in which $\Phi'X$ is the one-point compactification of
the space $U/G$, where
$$
U=\bigcup_{g\in G}
\{ a\in\map (K,X-\{ *\})\st \pi (a)=g\lambda\}.
$$
Now $U/G$ is homeomorphic to the orbit space
$$
\{ a\in\map (K,X-\{ *\})\st \pi (a)=\lambda\}/G_\lambda ,
$$
where $G_\lambda\leq G$ is the stabilizer of $\lambda$. 
Thus we can make the identification
$$
\Phi'X = {}_{G_\lambda}\Phi^{K_\lambda}_{\{ 0\}}X,
$$
where $K_\lambda$ is the set $\lambda$.

So we conclude from Lemma \ref{add}
that, if ${}_{G_\lambda}\Phi^{K_\lambda}_{\{ 0\}}$
is Lefschetz-polynomial of degree at most $\# K_\lambda$,
which is less than $k$, and ${}_G\Phi^K_{\LL^+}$ is Lefschetz-polynomial
of degree at most $k$, then ${}_G\Phi^K_\LL$
is Lefschetz-polynomial of degree at most $k$.

The proof is completed by induction, first on $k$, then on
the number of elements
$\# \PP (K)-\#\LL$ in the complement of $\LL$.
\end{proof}
\begin{proof}[Completing the proof of Theorem \ref{main}]
Take $G=\SS_k$ to be the group of all permutations of
$K=\{ 1,\ldots ,k\}$ and $\LL$ to be the set
of all partitions of $K$ into subsets with at most $l$ elements.
\end{proof}
\Cor{For any polynomial $p(t_1,\ldots ,t_k)\in\Qq [t_1,\ldots ,t_k]$
of degree $k$, there is a Lefschetz-polynomial functor $\Phi$
of degree $k$ and a positive integer $r$ such $L_\Phi =r\cdot p$.
}
\begin{proof}
Theorem \ref{main} with $l=1$ provides a Lefschetz-polynomial
functor $\Phi$ of degree 
$k$ with $L_\Phi(t_1,\ldots ,t_k)$ equal to the coefficient of
$q^k$ in $\prod_{m=1}^k (1+q^m)^{t_m}$, which is $t_k$ plus terms
involving $t_1,\ldots ,t_{k-1}$.
The result follows by induction on $k$ using wedge and smash products
and applying Lemmas \ref{add} and \ref{mult}.
\end{proof}

More generally, suppose that $Y$ is a pointed compact
$G$-ENR.
Define ${}_G\Phi^K_\LL (X; Y)$, regarded as a functor of $X$ for 
fixed $Y$, to be the quotient
$$
((\map (K,X)/(\{ a\st \pi (a)\notin\LL \}\cup \map_0(K,X)))\wedge Y)/G\, .
$$

When $\LL =\PP (K)$, following the proof of Theorem \ref{Gsymm}
we can compute the Lefschetz number as
$$
\frac{1}{\#G} \sum_g (-1)^i \tr ((\Bigwedge^kf^*)g^{-1}:
\tilde H^i(\Bigwedge^k X\wedge Y;\,\Qq )\to
 \tilde H^i(\Bigwedge^k X\wedge Y;\,\Qq )).
$$
In terms of $A_i=f^* : \tilde H^i(X;\,\Qq )\to \tilde H^i(X;\, \Qq )$,
the contribution from $g\in G$ is equal to the product
over the $g$-orbits $C$ of length $n_C$ of
$$
\sum_i (-1)^i \tr (A_i^{n_C})=\tilde L (f^{n_C})
=\sum_{m\, |\, n_C} m\tilde\DD_m(f)
$$
divided by $\# G$ and multiplied by the Lefschetz number
$\tilde L(g^{-1}: Y \to Y)$.
Hence the functor ${}_G\Phi_\LL^K(-;Y)$ is 
Lefschetz-polynomial in this case.

In general, we have a cofibre sequence
$$
{}_{G_\lambda}\Phi^{K_\lambda}_{\{ 0\}}(X; Y)
\to {}_G\Phi^K_{\LL^+} (X; Y) \to {}_G\Phi^K_\LL (X; Y)
$$
and deduce that ${}_G\Phi^K_\LL (-; Y)$ is Lefschetz-polynomial
for any $\LL$.

Specializing to symmetric powers, for a compact ENR $M$
and a pointed compact ENR $N$, we set
$$
S^k_l(M; N) ={}_{\SS_k}\Phi^K_{\LL}(X; \Lambda^kN),
\quad (k\geq 1,\, l\geq 0),
$$
where $K=\{1,\ldots ,k\}$,
$\SS_k$ acts on $K$ and the $k$-fold smash product $\Lambda^kN$
by permutation,
and $\LL \subseteq\PP (K)$ is the set of partitions into
subsets with at most $l$ elements.
We also write $S^k_\infty (M; N)=S^k_l(M; N)$ for $l\geq k$
and $S^0_l(M; N)=S^0$.
\Prop{\label{coeffic}
Let $N$ be a pointed compact ENR with reduced 
Euler characteristic $n\in \Zz$. Then for a self-map $f: M\to M$ of
a compact ENR
$$
\sum_{k\geq 0} \tilde  L(S^k_\infty (f) : S^k_\infty (M; N)\to
S^k_\infty (M; N))\, q^k =Z(f; q)^{-n}\in \Qq [[q]];
\leqno{{\rm (i)}}
$$
$$
\sum_{k\geq 0} \tilde L(S^k_1f : S^k_1(M; N)\to S^k_1(M; N))\, q^k
=\prod_{m\geq 1}(1+nq^m)^{\DD_m(f)}\, ;
\leqno{{\rm (ii)}}
$$
and, if $n\leq 0$,
$$
\sum_{k\geq 0} \tilde L(S^k_l(f) : S^k_l(M; N)\to
S^k_l(M; N))\, q^k =Z(f; q)^{-n}\in \Qq [[q]].
\leqno{{\rm (iii)}}
$$
for all $l\geq -n$.
}
\par\noindent
In particular, if $n=-1$, then the sum in (iii)
is equal to $Z(f; q)$ for any $l\geq 1$.
\begin{proof}
It is enough to check the identities when $M$ is
a finite set.
\end{proof} 

Suppose that $M$ is a connected, orientable, closed 
manifold of dimension $r$
and that $g :M \to M$ is a diffeomorphism. 
Write $\epsilon =+1$ if $g$ is orientation-preserving
and $\epsilon =-1$ if $g$ reverses orientation.

The configuration space $C^k(M)$ has the homotopy type of a 
finite complex.  We can, therefore, define the {\it Lefschetz trace}
of the induced map $C^kg : C^k(M)\to C^k(M)$ to be
$$
L^T(C^kg)=\sum_i (-1)^i\tr \{ (C^kg)^* : H^i(C^k(M);\,\Qq )
\to H^i(C^k(M);\,\Qq )\}.
$$
\Cor{In terms of the inverse $f=g^{-1}: M\to M$,
$$
\sum_{k\geq 0}  \epsilon^kL^T(C^kg)\, q^k=
\begin{cases}
Z(f; q)&\text{if $r$ is odd,}\\
Z(f; q^2)Z(f; q)^{-1}&\text{if $r$ is even.}
\end{cases}
$$
}
\begin{proof}
The configuration space $C^k(M)$ is the quotient
of the ordered configuration space $F_k(M)$
of $k$-tuples $(x_1,\ldots ,x_k)\in M^k$ with
$x_i\not=x_j$ for $i\not=j$ by the action of $\SS_k$.

Consider, in the notation of
Proposition \ref{coeffic} with $N=S^r$,  the
pointed space $S^k_1(M; S^r)$.
This is the one-point compactification of
the orientable $2kr$-dimensional manifold
$E_k=F^k(M)\times_{\SS_k} (\Rr^k\otimes\Rr^n)$. 
By Poincar\'e duality we have a non-singular pairing
(involving cohomology $H^*_c$ with compact supports)
$$
H^{kr+i}_c(E_k;\,\Qq )\otimes H^{kr-i}(E_k;\,\Qq )
\to H^{2kn}_c(E_k;\,\Qq )\iso\Qq
$$
such that $\langle f^*u,v\rangle = \epsilon^k\langle u,g^*v\rangle$,
for $u\in H^{kr+i}_c(E_k;\,\Qq )$, $v\in H^{kr-i}(E_k;\,\Qq )$.

The result follows, after making the identifications
$$
\text{
$H_c^{kr+i}(E_k;\, \Qq )=\tilde H^{kr+i}(S^k_1(M;\, S^r);\,\Qq )$
and $H^{kr-i}(E_k;\,\Qq )=H^{kr-i}(C^k(M);\,\Qq )$,}
$$
from Proposition \ref{coeffic} (ii) with $n=(-1)^r$
and $L(S^k_1f)=\epsilon^kL^T(C^kg)$.
\end{proof}
\Rem{If $M$ is a Riemannian manifold and $g$ is an isometry,
then $g$ induces, for any $\epsilon >0$,
a self-map of the compact subspace $C^k_\epsilon (M)$
of $C^k(M)$ consisting of those $k$-element subsets of $M$
in which any two elements are a distance $\geq\epsilon$
apart. For $\epsilon$ sufficiently small, the inclusion
$C_\epsilon^k(M)\into C^k(M)$ is a homotopy equivalence
and $L^T(C^kg)$ is equal to the Lefschetz number
of the restricted map $C^k_\epsilon (M)\to C^k_\epsilon (M)$.
} 
\Ex{Take $M$ to be the circle $S(\Cc )$ and
$g$ to be complex conjugation.
Thus, $r=1$, $\epsilon=-1$, $L(f)=2$ and $Z(f; q)=
(1-q)^2(1-q^2)^{-1}$.
We have $L^T(S^kg)=2$ if $k\geq 1$.
}
For $(1-q)^2(1-q^2)^{-1}=1-2q+2q^2-2q^3+\dots$.
(To check the Lefschetz numbers directly,
we can replace
$S^k(M)$ and $S^kg$ by $S(\Cc )$ with the conjugation
involution.)
\Sect{Powers}
Finally, let us specialize to the case in which the group $G$ is
trivial as studied in \cite{bary, hainaut}.

For a compact ENR $M$ and integers $k\geq 1$,
$l\geq 0$, let $F_l^k(M)\subseteq M^k$
be the subspace of points $(x_1,\ldots ,x_k)$ such that 
for all $x\in M$, $\#\{ i\st x_i=x\}\leq l$. 
It is an open subspace of $M^k$. In particular,
$F_1^k(M)$ is the ordered configuration space $F^k(M)$
of $M$, and $F_0^k(M)=\emptyset$. 
The quotient $T^k_l(M)=M^k/(M^k-F_l^k(M))$ is a pointed
compact ENR and may be identified with the one-point
compactification of $F_l^k(M)$. 
For $k=0$, we take $T^0_lM$ to be $S^0$.
A self-map $f$ of $M$ induces a pointed self-map $T^k_lf:
T^k_l M\to T^k_l M$.
\Thm{\label{sub}
{\rm (See \cite[Lemma 2.5]{hainaut}, \cite{BL}, 
\cite[Theorem 1.2]{bary} for the case that $f$ is the identity.)}
Let $f :M\to M$ be a self-map of a compact ENR.
Then, for $l\geq 0$,
$$
\sum_{k\geq 0} \tilde L(T^k_lf : T^k_lM\to T^k_lM )\, q_k
= (1+q_1+q_2+\cdots +q_l)^{L (f)}\in\Qq [[q]],
$$
where $q_i=q^i/i!$, so that $q_i\cdot q_j = \binom{i+j}{i,\, j}q_{i+j}$.
}
\begin{proof}
If $M$ is a finite set, the fixed points of $T^k_lf$, other than the 
basepoint $*$,
are the $k$-tuples $(x_i)$ in $\Fix (f)^k$ in which for
each $x\in\Fix (f)$ at most $l$ of the entries $x_i$ are equal
to $x$.
We see that
$$
\sum_k  \tilde L(T^k_lf)q_k =\prod_{x\in\Fix (f)} (1+q_1+\cdots +q_l).
$$

The general case follows from Proposition \ref{order} below
in the special 
case in which $K=\{ 1,\ldots ,k\}$ and $\LL$ is the family of partitions into blocks with cardinality at most $l$.
\end{proof}
\Prop{\label{order}
Consider a finite non-empty set $K$ with $\# K=k$
and a non-empty family $\LL\subseteq\PP (K)$ of
partitions closed under refinement.
Let $\Phi ={}_0\Phi^K_\LL$ be the functor defined in Section 3
for the case in which $G=0$ is the trivial group.
Then $\Phi$ is Lefschetz-polynomial of degree $k$
and $L_\Phi (t_1,\ldots ,t_k)=\ell^K_\LL (t_1)$,
where $\ell^K_\LL (t)\in\Zz [t]$ is the polynomial
$$
\ell^K_\LL (t)=\sum_{r=1}^k n_r(\LL )\, t(t-1)\cdots (t-r+1),
$$
expressed in terms of the number, $n_r(\LL )$,
of partitions in $\LL$ into $r$ disjoint non-empty subsets of $K$.
}
\begin{proof}
The functor is Lefschetz-polynomial by Lemma \ref{kl}.

It is enough, therefore, to calculate $\tilde L (\Phi f)$
when $X=M_+$ is a finite set. The Lefschetz number
counts the maps $a$ such that $a(i)\in\Fix (f)$ for all $i\in K$ and 
$\pi (a) \in \mathcal{L}$.
The number of such maps $a$ is equal to $\ell^K_\LL (\#\Fix (f))$.
\end{proof}
\Lem{{\rm (Compare \cite[2.5]{bary}.)}
Suppose that $K$ is a disjoint union $K=K_1\sqcup K_2$
of non-empty sets $K_i$ of size $k_i$
and $\LL$ has the form $\LL_1\times \LL_2$, where
$\LL_i\subseteq\PP (K_i)$. Then 
$$
\ell^K_\LL(t)=\sum_{r_1=1}^{k_1}\sum_{r_2=1}^{k_2}
n_{r_1}(\LL_1 )n_{r_2}(\LL_2) \, t(t-1)\cdots (t-(r_1+r_2)+1)\, .
$$
}
\begin{proof}
This just says that $n_r(\LL )=\sum_{r_1+r_2=r} n_{r_1}(\LL_1)
n_{r_2}(\LL_2)$.
\end{proof}
\Ex{{\rm (Compare \cite[Example 1.6]{bary}.)}
Suppose that $K=K_1\sqcup K_2$ where $\# K_1=j=\# K_2$ and $\LL$ is 
the set of all refinements of the partition  $\{K_1,\, K_2\}$. 
Then
$$
\ell^{K}_\LL (t)=
\sum_{r_1,r_2=1}^j S(r_1,j)S(r_2,j)\, t(t-1)\cdots (t-(r_1+r_2-1)),
$$
where $S(r,j)$ is the Stirling number counting the number of
partitions of a $j$-element set into $r$ blocks.
For small positive values of $t$ we have:
$\ell^{K}_\LL (1)=0$;
$\ell^{K}_\LL (2)=2$ if $j\geq 1$;
$\ell^{K}_\LL (3)=6(2^j-1)$ if $j\geq 2$.
}
\begin{appendix}
\Sect{Cohomology of the orbit space}
\Prop{Let $G$ be a finite group and $M$ a compact $G$-ENR.
Then the projection $\pi : M\to M/G$ induces an
isomorphism $H^*(M/G;\, \Qq )\to H^*(M;\, \Qq )^G$.
}
\begin{proof}
We shall show more precisely that, for any prime $p$, including $p=0$,
which does not divide $\# G$,
$$
\pi^* : H^i(M/G;\, \Zz )_{(p)} \to H^i(M;\,\Zz )_{(p)}^G
\quad (i\in\Zz )
$$
is an isomorphism.

Consider first the case $i=0$. The cohomology $H^0(M;\,\Zz )$ is 
the abelian group of continuous maps $M\to \Zz$.
And $H^0(M/G;\,\Zz )$ is the group of continuous maps
$M/G\to\Zz$. So clearly $H^0(M/G;\,\Zz )=H^0(M;\,\Zz )^G$.

For $i\geq 1$, we use the Dold-Thom theorem. As in 
\cite[II, Section 15]{FHT}, we write $A(N)$ for the free commutative monoid on a connected pointed compact ENR $N$ with the basepoint 
as zero element.
It is filtered by subspaces $A^k(N)$ consisting of the sums of at most 
$k$ elements of $N$, topologized as a quotient of $S^kN$.
(In the notation employed in the proof of Theorem \ref{Gsymm},
$A^k(N)$ is the quotient of $\map (K,N)/\SS_k$ by the subspace
$\map_0(K,N)/\SS_k$, where $K=\{ 1,\ldots ,k\}$.)
The Dold-Thom theorem identifies $H^i(M;\,\Zz )$, for $i\geq 1$,
with the monoid of homotopy classes of maps $a:M\to A(S^i)$.
This allows us to construct a {\it transfer} homomorphism
$$
\pi_* : H^i(M;\, \Zz )\to H^i(M/G;\,\Zz )
$$
sending the class of a map $a$ to the class of the map
$\pi_*(a): M/G\to A(S^i)$:
$[x]=\pi (x)\mapsto \sum_{g\in G} a(gx)$, $x\in M$.
If $a$ maps into $A^k(S^i)$, $\pi_*(a)$ maps into $A^{k\# G}(S^i)$.
From the definition it is clear that
$$
\text{$\pi^*\pi_*(a) =\sum_{g\in G} ga$ for $a\in H^i(M;\,\Zz )$
and $\pi_*\pi^*(b)= \# G \cdot b$ for $b\in H^i(M/G;\,\Zz )$.}
$$
Notice that $\pi^*\pi_*(a)$ lies in $H^i(M;\,\Zz )^G$, and, 
if $a\in H^i(M;\,\Zz )^G$, is equal to $\# G\cdot a$.

At the prime $(p)$ we can invert $\# G$, and $(\# G)^{-1}\pi_*$
is an inverse for $\pi^*$ on the invariant subgroup 
$H^i(M;\,\Zz )_{(p)}^G$.
\end{proof}
\end{appendix}

\end{document}